\documentstyle {article}

\title{Morita equivalence of multidimensional noncommutative
tori.}
\author{Marc A. Rieffel and Albert Schwarz
\thanks{M. A. R. supported in part by 
NSF grant DMS 96-13833; A. S. supported in part by NSF grant DMS
95-00704; both authors supported in part by NSF Grant No. PHY94-07194 
during visits to the Institute for Theoretical Physics, 
Santa Barbara. A.S. acknowledges also the hospitality of MIT, IAS and
Rutgers University }
\\Department of Mathematics, UC Berkeley, Berkeley, 
CA 94720-3840 
\\Department of Mathematics, UC Davis, Davis, CA 95616}
\date{June 18, 1998}
\begin {document}
\maketitle
\smallskip

\begin{abstract}
One can describe an $n$-dimensional noncommutative torus by means of an
antisymmetric $n\times n$ matrix $\theta$. We construct an action of the group
$SO(n,n|{\bf Z})$ on the space of $n\times n$ antisymmetric matrices and
show that, generically, matrices belonging to the same orbit of this group 
give Morita equivalent tori. Some applications to physics are sketched. 
\end{abstract}
 
By definition \cite{R5}, 
an $n$-dimensional noncommutative torus is an associative
algebra with involution having unitary generators $U_1,...,U_n$ obeying
the relations 
\begin {equation}
U_kU_j=e(\theta _{kj})U_jU_k  ,
\end {equation}
where $e(t) = e^{2\pi it}$ and $\theta$ is an antisymmetric matrix. 
The same name is used for different completions of this algebra. In
particular, we can consider the noncommutative torus as a $C^{\star}$-algebra
$A_{\theta}$ (the universal $C^{\star}$-algebra generated by $n$ unitary
operators satisfying (1) ). Noncommutative tori are important in many
problems of mathematics and physics. It was shown recently that they are
essential in consideration of compactifications of M(atrix) theory
(\cite{CDS}; for further development see \cite{T}).
The results of the present paper also have application 
to physics. 

If two algebras ${\cal A}$ and $\hat {{\cal A}}$ 
are Morita equivalent (see the definition below), then for every ${\cal
A}$-module $R$ 
one can construct an $\hat {{\cal A}}$-module $\hat {R}$ 
in such a way that the correspondence $R\rightarrow \hat {R}$ 
is an equivalence of categories of ${\cal A}$-modules and 
$\hat {{\cal A}}$-modules. 

It was shown in \cite{CDS} that in the consideration of toroidal 
compactification of M(atrix) theory one can identify states 
with connections on projective ${\cal A}$-modules, where 
${\cal A}$ is a noncommutative torus. 
 
One can verify that compactifications on Morita equivalent 
tori are in some sense physically equivalent; this fact is 
conjectured in \cite{CDS} and proved in \cite{S}.

In particular, it is shown  in \cite{S} that there is a correspondence
between 
BPS states on the module $R$ with BPS states on the module 
$\hat {R}$. (BPS states having 
maximal supersymmetry correspond to connections with 
constant curvature.  The results of \cite{S} can be applied also to the
case of BPS states 
having less supersymmetries .)

It was mentioned in \cite{CDS} that in the case of two-dimensional 
tori Morita equivalence is related to $T$-duality in string 
theory. This remark leads to a conjecture that the group 
$SO(n,n |{\bf Z})$ that appears in the present paper is related 
to a corresponding group  in $T$-duality.

It is clear that the relations (1) are unchanged if each $U_j$ is 
changed by some phase factor. One can use this fact to arrange in various
ways \cite{R4} 
that there is a bicharacter $\gamma$ from ${\bf Z}^n$ to the circle 
group ${\bf T}$, and unitary elements $U_x$ of $A_\theta$ for 
each $x\in {\bf Z}^n$ such that
\begin{equation}
U_x U_y = \gamma(x, y)U_{x+y}    ,
\end{equation}
and such that if $e_j$ denotes the $j^{\rm th}$ standard basis vector 
in ${\bf Z}^n$ and if $U_j = U_{e_j}$, then the relations (1) are satisfied. 
One useful choice for $\gamma$ is
$$
\gamma(x,y) = e((x\cdot \theta y)/2)  .
$$
In any event we will have 
$$
\gamma(x,y){\bar \gamma}(y,x) = \rho(x,y)
$$ where $\rho$ is the skew-bicharacter defined by
$$
\rho(x,y) = \rho_\theta(x,y)= e(x\cdot\theta y).
$$
We will use $\rho$ below.

One can consider also a smooth version
of a noncommutative torus \cite{R4, R5}, which is 
the algebra $A_{\theta}^{\infty}$ consisting of
formal series 
$$\sum c_{j_1...j_n}U_1^{j_1}...U_n^{j_n}$$
where the collection of coefficients belongs to the Schwartz 
space ${\cal S}({\bf Z}^n)$. One can characterize
$A_{\theta}^{\infty}$ as the subalgebra of ${A}_{\theta}$ consisting of
vectors which are smooth with respect to the natural action of the
$n$-dimensional commutative Lie algebra on ${A}_{\theta}$; this
action is given by *-derivations $\delta_j$ defined by 
the formula $\delta_j U_j=iU_j,\  \delta_j U_k=0$ if
$j\not= k.$ We will consider the  noncommutative torus as a
$C^{\star}$-algebra ${A}_{\theta}$; however, our results remain
correct if $A_{\theta}$ is replaced with $A_{\theta}^{\infty}$. 

The notion of (strong) Morita equivalence of $C^{\star}$-algebra was
introduced and analyzed in \cite{R1, R3}. 
We will use the following constructive
definition of Morita equivalence for the case of unital $C^{\star}$-algebras. 
Let us consider a finite projective right
module ${\cal E}$ over a $C^{\star}$-algebra ${A}$ (i.e. a module that
can be considered as a direct summand of a finite-dimensional free module
over ${A}$.) The algebra $End_{{A}}{\cal E}$ of endomorphisms of
${\cal E}$ has a canonical structure as a $C^{\star}$-algebra. We say that  a
$C^{\star}$-algebra  ${A}^{\prime}$ is (strongly) Morita equivalent
to ${A}$ if it is isomorphic to $End _{{A}}{\cal E}$ for
some finite projective module ${\cal E}$. (This is sufficient for our
present purposes, but more generally one must require that the module is full
\cite{R3}, which is automatic here.)

Morita equivalent  $C^{\star}$-algebras share many important properties.
(They have equivalent categories of modules, isomorphic $K$-groups, cyclic
homology etc.) 

Let us start our development 
with the definition of an action of the group $O(n,n|{\bf R})$
on the space ${\cal T}_n$  of real antisymmetric matrices. We consider
$O(n,n|{\bf R})$ 
as a group of linear transformations of the space ${\bf R}^{2n}$
preserving the quadratic form $x_1x_{n+1}+x_2x_{n+2}+...+x_nx_{2n}$.
It is convenient to denote coordinates in ${\bf R}^{2n}$ as
$(a^1,...,a^n,b_1,...,b_n)$. Then the quadratic form on ${\bf R}^{2n}$ can
be written as $a^ib_i$. Accordingly, we will often write the
elements of $O(n,n|{\bf R})$ in $2\times 2$ block form
$$
g = \left( \begin{array}{cc}
A & B\\
C & D
\end{array} \right).
$$
Here the blocks $A,B,C,D$ are $n\times n$ matrices which 
satisfy $A^tC+C^tA=0=B^tD+D^tB,\  \  A^tD+C^tB=1$, where $^t$ 
denotes transpose. The action of
$O(n,n |{\bf R})$ on the space ${\cal T}_n$ of antisymmetric matrices is
defined by the formula 
\begin{equation}
\theta ^{\prime}= (A\theta + B)(C\theta +D)^{-1}  .
\end{equation}
We emphasize that this action is defined only on the 
subset of ${\cal T}_n$ where $(C\theta + D)$ is invertible. This
subset can be empty, e.g. for $\sigma_1$ as defined before the lemma
below. But we will see shortly that this set is dense for every
element in the subgroup $SO(n,n|{\bf Z})$ of $O(n,n|{\bf R})$
consisting of matrices with integer entries and determinant $+1$, which is
the subgroup we need to use anyway for other reasons.
 
In order to give an elegant formulation of our results, we let
$$
{\cal T}_n^0 = \{\theta \in {\cal T}_n:g\theta {\rm \ is\  defined\ 
for\  all \ } g \in SO(n,n|{\bf Z})\}    .
$$
A simple calculation shows that ${\cal T}_n^0$ is carried into itself by
the action of  
$SO(n,n|{\bf Z})$. But it is not so obvious how big ${\cal T}_n^0$ is.
However, the considerations which we must go through anyway in order
to prove our main theorem will give us a way of seeing that ${\cal T}_n^0$
already 
is ``big'' in a suitable sense, and that in particular it is dense in ${\cal
T}_n$.
\medskip

{\bf Theorem.} For  $\theta \in {\cal T}_n^0$ and
$g\in SO(n,n|{\bf Z})$ the  noncommutative torus corresponding 
to $g\theta $ is
Morita equivalent  to the  noncommutative torus  corresponding to $\theta$.
This remains true for the smooth versions of the non-commutative tori.
\medskip

We do not know whether the converse holds at the level of the smooth
algebras, that is, whether if 
two $\theta$'s give Morita equivalent smooth noncommutative tori 
then these $\theta$'s must be in the same orbit for the 
action of $SO(n,n|{\bf Z})$. But this converse is false at the C*-algebra
level, as we discuss after the end of the proof of the theorem. 
However, when $n = 2$ the converse does
hold at the C*-algebra level 
\cite{R2}, and the appearance here of $SO(n,n|{\bf Z})$ is
a generalization of the appearance of $SL(2, {\bf Z})$ discovered
in \cite{R2}. In the case $n>2$ the converse statement is proved in
\cite{S}
for a modified definition of Morita equivalence.
We will comment on what one can say for $\theta$'s which are not in
${\cal T}_n^0$ after the end of the proof of the theorem.

For the proof of our this theorem
we now define a suitable set of generators for $SO(n,n|{\bf Z})$.

For every matrix $R \in GL(n|{\bf Z})$ 
we can define the 
transformation $\rho (R)\in SO(n,n|{\bf Z})$ by the formula 
$$
\tilde {a}^i=R_j^ia^j,\ \ \ \ \ \tilde {b}_i=(R^{-1})_i^jb_j  .
$$
If $N$ is an anti-symmetric $n\times n$ matrix 
with integer entries $n_{ij}$ we can
construct an element $\nu (N)$ of $SO(n,n|{\bf Z})$ as a transformation 
$$
\tilde {a}^i=a^i +n^{ij}b_j,\ \ \ \ \   \tilde{b}_i=b_i.
$$
Finally, for every integer $k$ we define an 
element $\sigma _k\in O(n,n|{\bf Z})$
by the formula 
$$\tilde {a}^i=b_i \ \  {\rm for}\ \ 1\leq i\leq k,\ \ \ \tilde {a}^i=a^i \
\  {\rm for} \ \ k<i\leq n,$$
$$\tilde {b}_i=a^i \ \  {\rm for}\ \ 1\leq i\leq k,\ \ \ \tilde {b}_i=b_i \
\  {\rm for} \ \ k<i\leq n.$$
Notice that $\sigma _k \in SO(n,n|{\bf Z})$ only when $k$ is even. 
When $n$ is odd no element of ${\cal T}_n$ can be invertible, and so
$\sigma _n$ will have no element of  ${\cal T}_n$ in its domain.

\medskip

{\bf Lemma.} The elements $\rho (R),\ \  \nu (N)$, and the single
element $\sigma _{2}$, together
generate the group $SO(n,n|{\bf Z})$.
\medskip

Proof. Let $H$ denote the group generated by these elements. 
Clearly $H \subseteq SO(n, n|{\bf Z})$. For any anti-symmetric $n\times n$ 
matrix $N$ let $\mu(N)$ denote the transformation
$$
        \tilde {a}^i=a^i,\ \ \ \ \   \tilde{b}_i=b_i+n_{ij}a^j.    
$$
Now suppose that $N$ has its only non-zero (integer) entries in the 
top-left $2\times 2$ block. Then it is easily verified that when 
$\nu(N)$ is conjugated by $\sigma_2$ one obtains exactly 
$\mu(N)$. Thus $\mu(N) \in H$. If we now conjugate $\mu(N)$ by 
$\rho(R)$ as $R$ ranges over all permutation matrices, we obtain 
all the $\mu(N)$ for all $N$ which have only 2 non-zero entries. 
But every anti-symmetric matrix is a sum of such $N$'s, and the 
map $N \mapsto \mu(N)$ is a homomorphism from the additive 
group of integer-valued anti-symmetric matrices into $SO(n,n|{\bf 
Z})$. Thus $H$ contains all $\mu(N)$ for all such matrices $N$.

Let $EO(n,n|{\bf Z})$ denote the subgroup of $SO(n,n|{\bf Z})$ 
generated by all the $\mu(N)$'s and $\nu(N)$'s, together with 
those $\rho(R)$ for which $R$ is an elementary matrix (i.e. is $I_n$ 
plus only one off-diagonal integer non-zero entry). Then 
$EO(n,n|{\bf Z})  \subseteq  H$.  Notice that $O(1,1|{\bf Z})$ is a 4-element 
group, with subgroup $SO(1,1|{\bf Z})$ of order 2. When this 
observation is used in theorem 5.5.3 of \cite{HO}, one finds that 
$EO(n,n|{\bf Z})$ is a normal subgroup of $SO(n,n|{\bf Z})$ of index 
at most 2, and that $EO(n,n|{\bf Z})$ together with $\rho(U)\oplus 
I_{n-2}$ generates $SO(n,n|{\bf Z})$, where $U$ is the non-identity 
element of $SO(1,1|{\bf Z})$. (See the remark near the bottom of 
page 232 of \cite{HO}.) Since this element is in $H$, it follows that $H$ 
= $SO(n,n|{\bf Z})$ as desired. (From theorem 7.2.23 of \cite{HO} 
one sees that $EO(n,n|{\bf Z})$ is not itself all of $SO(n,n|{\bf Z})$.)     
\medskip

We now show that ${\cal T}_n^0$ is ``big''.
Let $G$ denote the subgroup of $SO(n,n|{\bf Z})$ 
generated by the $\rho(R)$'s and the $\nu(N)$'s, so that it consists 
exactly of the $g \in SO(n,n|{\bf Z})$ for which $C = 0$. It is clear 
that $G$ acts on all of ${\cal T}_n$. Note that the action of $\sigma = 
\sigma_2$ is defined on $\theta \in {\cal T}_n$ exactly if the top-left 
$2 \times 2$ block of $\theta$ is non-zero, and so invertible. Let 
$U$ denote the set of such $\theta$'s. Then $U$ is an open dense 
subset of ${\cal T}_n$, which is carried onto itself by $\sigma$ and
contains ${\cal T}_n^0$. Let 
$$
V_1 = \bigcap \{gU: g \in SO(n,n|{\bf Z})\}          .
$$
Since $SO(n,n|{\bf Z})$ is countable, Baire's theorem tells us that 
$V_1$ is a dense subset of ${\cal T}_n$ of the ``second category" 
(that is, it is not the countable union of nowhere-dense sets). 
Clearly $V_1$ is $G$-invariant and 
contains ${\cal T}_n^0$. Let $U_1 = \sigma(V_1)$. Then let
$$
V_2 = \bigcap \{gU_1: g \in SO(n,n|{\bf Z})\}\cap V_1       ,
$$
and let $U_2 = \sigma(V_2)$. We continue inductively to define 
$V_k$ for all $k > 0$. Each $V_k$ is $G$-invariant and
contains ${\cal T}_n^0$. Set ${\cal T}_n^1 = 
\bigcap^{\infty} V_k$. Again ${\cal T}_n^1$ is dense of the second 
category, so is ``big". Clearly ${\cal T}_n^1$ is $G$-invariant
and
contains ${\cal T}_n^0$. But if 
$\theta \in {\cal T}_n^1$, then $\theta \in V_{k+1} \subseteq U_k$ 
for each $k$, so that $\sigma(\theta) \in V_k$. Thus $\sigma(\theta) 
\in {\cal T}_n^1$, and so $\sigma$ carries ${\cal T}_n^1$ into itself. In 
view of the Lemma, $SO(n,n|{\bf Z})$ has a fully-defined action on 
${\cal T}_n^1$, carrying ${\cal T}_n^1$ into itself. Since ${\cal T}_n^1$
contains  ${\cal T}_n^0$, it follows from
the definition of ${\cal T}_n^0$ that ${\cal T}_n^1 = {\cal T}_n^0$.
Thus ${\cal T}_n^0$ is ``big'' as desired.

\medskip

It is clearly sufficient to prove the theorem just for the $g$'s in
the generating set for $SO(n,n|{\bf Z})$ given in the Lemma above.
To analyze the case $g=\rho (R)$ we should check that  the
noncommutative torus determined by the matrix $\theta ^{\prime}=R\theta
R^t$  is Morita equivalent  to the  torus  corresponding to the matrix
$\theta$. One can verify that these tori are in fact isomorphic. This is 
seen by using the description given above of noncommutative tori 
in terms of bicharacters. The matrix $R\in SL(n,{\bf Z})$ generates an
automorphism of ${\bf Z}^n$. The isomorphism between ${A}_{\theta
^{\prime}}$ and ${A}_{\theta}$, and between the smooth versions, 
follows from the relation 
\begin {equation}
U_{R(x)}U_{R(y)}=e((x\cdot {\theta}^{\prime }y)/2)U_{R(x+y)}.
\end {equation}

Next, it is obvious that upon replacing $\theta _{ij}$ with $\theta
_{ij}^{\prime}=\theta _{ij}+n_{ij}$, where the $n_{ij}$ are integers, we do
not  change the commutation relations (1) at all. This means that
$\theta$
and $\theta ^{\prime}=\nu (N)\theta$  correspond to the  same
noncommutative torus.

It remains to prove that ${A}_{\theta}$ is Morita equivalent  to
${A}_{\theta ^{\prime}}$ when $\theta ^{\prime}=\sigma _{2}\theta
$. We will in fact prove this for any $\sigma _{2p}$, and any $\theta$
in its domain.
\medskip

{\bf Proposition.} Let $\theta \in {\cal T}_n$, and assume that the 
top left $2p\times 2p$ block of $\theta$ (which we will denote 
by $\theta_{11}$) is invertible. Write $\theta$ in block form as

$$
 \theta = \left( \begin{array}{cc}
\theta_{11} & \theta_{12}\\
\theta_{21}  & \theta_{22}
\end{array} \right) .
$$
Let

$$
\theta'  =  \sigma_{2p}(\theta)  =  
\left( \begin{array}{cc}
\theta_{11}^{-1}  &  -\theta_{11}^{-1}\theta_{12}\\
\theta_{21}\theta_{11}^{-1}  &  
\theta_{22} - \theta_{21}\theta_{11}^{-1}\theta_{12}
\end{array} \right)         .
$$
Then the non-commutative torus for $\theta'$ is Morita equivalent to that 
for $\theta$. This remains true at the level of the smooth algebras.
\medskip

Proof. For the proof we use the construction of projective modules given in 
\cite{R4}. These modules are constructed at the smooth level, and so our
proof here works at both the smooth level and the C*-algebra level. 
To make the argument more transparant, we use ordering 
conventions which 
are especially adapted to our present situation. These conventions are 
slightly different from those used in \cite{R4}, 
but they only affect matters of 
orientation. Since we will not use the apparatus of connections and Chern 
character developed in \cite{R4}, this difference will have no effect on our 
present considerations.

Let 
$$
J_o = \left( \begin{array}{cc}
0 & I_p\\
-I_p &  0
\end{array} \right)     .
$$
Then $-\theta_{11}$ is similar to  $J_o$, and so we can choose an invertible 
matrix, $T_{11}$, such that  $T_{11}^t J_o T_{11} = -\theta_{11}$. 
Set $T_{13} = 
\theta_{12}^t$. Set $q = n-2p$ and let $T_{32}$ be any $q\times q$ 
matrix such that  
$\theta_{22} = T_{32}^t - T_{32}$. For example, $T_{32}$ can be the part of 
$\theta_{22}$ above the main diagonal, with $0$'s below, or alternatively it 
can be $-\theta_{22}/2$.

Set
$$
T = \left( \begin{array}{cc}
T_{11}  &  0\\
0  &  I_q\\
T_{31}  &  T_{32}
\end{array} \right)            ,
$$
a matrix of size $(n+q)\times n$. Also set
$$
J =  \left( \begin{array}{ccc}
J_o  &  0  &  0\\
0  &   0  &  I_q\\
0  &  -I_q  &  0
\end{array} \right)         ,
$$
a square matrix of size $n+q$, which is for us a convenient way to write 
the standard symplectic matrix. A routine calculation shows that $T^t JT = 
-\theta$. (The minus sign is included because, as required by the 
definition of Morita equivalence, we will construct a right module. If one 
arranges for $+\theta$ here, then there will be unpleasant signs in the final 
formula.) Following the notation in \cite{R4}, we let 
$$
{\tilde T} = \left( \begin{array}{cc}
T_{11}  &  0\\
0  &  I_q
\end{array} \right)         ,
$$
and we note that ${\tilde T}$ is invertible. Note also that 
as a linear transformation, $T$ carries ${\bf Z}^{2p}\times {\bf Z}^q$ into 
${\bf R}^{2p}\times {\bf Z}^q\times {\bf R}^q$. Thus $T$ satisfies the 
conditions in definition 4.1 of \cite{R4} for being an ``embedding" map. This 
means, as is anyway clear, that when we view $T$ as a homomorphism from 
${\bf Z}^n$ into  $G = {\bf R}^{2p}\times {\bf Z}^q\times {\bf T}^q$, 
then the 
range of $T$ is a lattice in $G$, which we will denote by $D$. 
In what follows,
our notation will not distinguish between elements of ${\bf R}^{2p}\times 
{\bf Z}^q\times {\bf R}^q$  and their images in $G$.  

We view $G$ as 
$M\times \hat M$ with $M = {\bf R}^p\times {\bf Z}^q$, where $\hat {}$ 
denotes ``dual group". From this decomposition $G$ carries a canonical 
cocycle, $\beta$, called the ``Heisenberg cocycle'' in \cite{R4}, and
defined by
$$
\beta ((m,s), (n,t)) = \langle m, t\rangle   ,
$$
where here $\langle\  , \  \rangle$ denotes 
the pairing between $M$ and $\hat M$. The corresponding 
skew cocycle, $\rho$, is defined by
$$
\rho ((m,s), (n,t)) = \langle m, t\rangle\overline {\langle n,s\rangle}.
$$
It is clear that (after rearranging the order of the 
factors) we have exactly $\rho(x, y) = e(x\cdot Jy)$ for the $J$ 
defined above. 
Because $T^tJT = -\theta$, the restriction of $\rho$ to 
$D$ and so to ${\bf Z}^n$ is exactly $(x,y) \mapsto e(-x\cdot \theta y)$.
Except for a complex conjugation, this is exactly the $\rho$ defined much 
earlier on ${\bf Z}^n$. 

The restriction of $\bar \beta$ to $D$ will then be a cocycle whose 
corresponding skew cocycle is exactly the restriction of $\rho$ to $D$,
and so to ${\bf Z}^n$. Thus $\bar \beta$ is an instance of the cocycle
$\gamma$ of equation (2) above. Consequently, when we
let $A = {\cal S}({\bf Z}^n, \bar \beta)$, 
the space of Schwartz 
functions on ${\bf Z}^n$ with convolution twisted by the cocycle
$\bar \beta$, this algebra, or rather its $C^{\star}$-algebra completion, 
is (isomorphic to) our non-commutative torus $A_{\theta}$. 

Because of the minus sign in $e(-x\cdot \theta y)$, $A$ has a natural
{\it right} action on the space ${\cal S}(M)$ of Schwartz functions
on $M$, defined by a twisted convolution indicated shortly after the
proof of proposition 2.9 of \cite{R4} (except that there $D^\perp$ plays the
role of our $D$).
We do not need the formula for this action because what we really need is a 
description of the endomorphism algebra of this right $A$-module, since it is 
this endomorphism algebra which is Morita equivalent to $A$ via this 
module. But by applying proposition 3.2 of \cite{R4} we see that this 
endomorphism algebra is exactly ${\cal S}(D^\perp , \beta)$ where
$$
D ^\perp = \{w\in G: \rho(w, z) = 1 {\rm \ \ for \  all \ \ } z \in D\}   ,
$$ 
and then we restrict $\beta$ to $D ^\perp$. This restricted $\beta$ will,
of course, have the restricted $\rho$ as its corresponding skew cocycle.

Thus we need to determine $D ^\perp$. More specifically, we need to describe 
$D ^\perp$ as the image of ${\bf Z}^n$ under some embedding map, so that we 
can calculate the anti-symmetric matrix $\theta'$ which on ${\bf Z}^n$ 
gives the cocycle corresponding to the restriction of $\rho$ to $D^\perp$. 
To obtain the answer 
specified in the statement of the theorem, we need to choose this embedding 
map quite carefully. We proceed as follows.

Given $x\in G$, it will be in $D ^\perp $ exactly 
if $x\cdot JTz \in {\bf Z}$ for 
all $z\in {\bf Z}^n$, that is, exactly if $T^tJx \in {\bf Z}^n$. 
We want a natural 
isomorphism from ${\bf Z}^n$ to this $D ^\perp$. To see how to obtain one, let
$$
\bar T = \left( \begin{array}{ccc}
T_{11}  &  0  &  0\\
0  &  I_q  &  0\\
T_{31}  &  T_{32}  &  I_q
\end{array} \right)          ,
$$
a square matrix of size $n+q$ which is clearly invertible. 
It is easily checked 
that $T^tJx \in {\bf Z}^n$ exactly if $\bar T^tJx \in {\bf Z}^{n+q}$. 
Since $\bar 
T$ is invertible, as is $J$, this means that
$$
                  D ^\perp = (\bar T^tJ)^{-1}({\bf Z}^{n+q})  ,
$$
viewed in $G$. But a routine calculation shows that
$$
(\bar T^tJ)^{-1} = \left( \begin{array}{ccc}
-J_o(T_{11}^t)^{-1}  &  0  &   J_o(T_{11}^t)^{-1}T_{31}^t \\
0  &  0  &  -I\\
0  &  I  &  -T_{23}^t
\end{array} \right)     .
$$
It is clear from this that
$$
(\bar T^tJ)^{-1}(0\times {\bf Z}^q\times 0) = 0\times 0\times {\bf Z}^q    ,
$$
which is $0$ in $G = {\bf R}^{2p}\times {\bf Z}^q\times {\bf T}^q$.
Thus we can omit the second column of $(\bar T^tJ)^{-1}$ . This gives us, 
after we omit an inessential sign,
$$
S =  \left( \begin{array}{cc}
J_o(T_{11}^t)^{-1}  &   -J_o(T_{11}^t)^{-1}T_{31}^t \\
0  &  I\\
0  &  T_{23}^t
\end{array} \right)     .    
$$
This is the desired embedding map giving an isomorphism 
from ${\bf Z}^n$ onto $D^\perp$.

To conclude, we must calculate the matrix for the cocycle on ${\bf Z}^n$ 
coming from $\rho
$ via $S$. But a routine calculation shows 
that $S^tJS =  \sigma_{2p}(\theta)$, whose matrix is given 
in the statement of the proposition.  

We mention that F. Boca has pointed out to us that the case of the
above proposition for $\sigma_2$ is implicit in calculations occurring
shortly before lemma 2.1 of \cite{B}, upon specialization of parameters.
See also the remarks following lemma 2.6 of \cite{B}.

Let us make some comments about the proof of our theorem, as it
applies to $\theta$'s which are not in ${\cal T}_n^0$. We
notice that each of our generators for $SO(n,n|{\bf Z})$ is  defined on 
an open dense
subset of ${\cal T}_n$. 
Let us suppose that $g \in SO(n,n|{\bf Z})$ is 
represented as $g_1,...,g_n$ where $g_i$ are from among these
generators. It is clear that we can find a dense open subset ${\cal
T}^g$ of ${\cal T}_n$ such that $g_i$ is well defined on $g_{i+1}...g_n{\cal 
T}^g$ for each $i$. 
It follows from our consideration that for $\theta \in {\cal T}^g$
the torus ${\cal A}_{g\theta}$ is Morita equivalent to ${\cal
A}_{\theta}$.

Let us suppose
that $\theta$ is in ${\cal T}_n$ but not in ${\cal T}^g$, 
for a given $g \in SO(n,n|{\bf Z})$, and that
$C \theta + D$ is invertible, so that $g\theta$ is defined. One
can conjecture  that in this case the noncommutative torus for $g \theta$
is also Morita
equivalent to that for $\theta$. We were not able to prove this
conjecture.

We now discuss the difference between the smooth and C*-algebraic
aspects. It is not difficult to see that a Morita equivalence between
smooth non-commutative tori gives a Morita equivalence between the
corresponding C*-algebras. But it is possible for two isomorphic
C*-algebraic non-commutative tori to have non-isomorphic smooth algebras.
Put another way, a given C*-algebraic non-commutative torus can have
inequivalent smooth structures. That this can happen follows from a deep 
investigation of the structure of 3-dimensional non-commutative tori
by G. A. Elliott and Q. Lin, culminating in \cite{L}, combined with a
deep investigation a decade earlier of isomorphisms of the smooth algebras,
culminating in \cite{BC}. In \cite{L} it is shown that the C*-algebraic
non-commutative 3-tori are entirely classified up to isomorphism by
their ordered $K_{0}$-group with distinguished order-unit. The Morita
equivalence classes are completely classified by the ordered $K_{0}$-group,
ignoring the order-unit. Furthermore, the Morita equivalence class is 
entirely determined in terms of the range of the canonical trace on the
$K_{0}$-group (corollary 2 of \cite{L}). 

The following
case is relevant to our present considerations. Choose  
$\theta \in {\cal T}_3^0$ such that the 7 numbers  $1, \theta_{12}, 
\theta_{13}, \theta_{23}$, together with all products of any two of
them, are linearly independent over the rational numbers. Since ${\cal T}_3^0$
is of second category, this can be done. Let $\phi$ be $\theta$ except with
$\theta_{12}$ replaced by $-\theta_{12}$. By the results mentioned in the
previous paragraph, the C*-algebras $A_\theta$ and $A_\phi$ are 
isomorphic, and so
are Morita equivalent.

However, we now show that $\phi$ is not in the orbit of $\theta$ for the
action of $SO(3,3|{\bf Z})$. If it were in the orbit, there would be
a $g \in SO(3,3|{\bf Z})$, given in block form as earlier, such that
  $$A\theta + B = \phi(C\theta + D) = \phi C \theta + \phi D.$$
By the linear independence assumed above, $B=0$ and $A\theta = \phi D$.
Because $B=0$, we must have $A^t D = I$. Consequently, $A\theta A^t = \phi$.
View $\theta$ and $\phi$ as defining antisymmetric bilinear forms on
${\bf Z}^3$. The preceeding equation says that $A$ carries the form for
$\theta$ to that for $\phi$. Now view these forms as defining linear
functionals, again denoted by $\theta$ and $\phi$, on 
${\bf Z}^3 \wedge {\bf Z}^3$. Then $\phi = \theta \circ (A\wedge A)$ 
on ${\bf Z}^3 \wedge {\bf Z}^3$. In particular, for the standard basis
of ${\bf Z}^3$ we have
  $$\epsilon_{jk}\theta_{jk} = \phi(e_j\wedge e_k) = 
                         \theta((A\wedge A)(e_j\wedge e_k)),$$
where $\epsilon_{12} = -1$ while $\epsilon_{jk} = +1$ for
$(j,k) = (1,3), (2,3)$. By the linear independence assumed above
for $\theta$ it follows that
  $$(A\wedge A)(e_j\wedge e_k) = \epsilon_{jk} e_j\wedge e_k .$$
Consequently, $\det(A\wedge A) = -1$. But because we are in dimension
3, $\det(A\wedge A) = (\det (A))^2 \geq 0$, which is a contradiction.

Along the same lines it can be seen as in \cite{BC} that the smooth
algebras for $\theta$ and $\phi$ are not isomorphic.

\medskip

{\bf Appendix.}

  One can define the action of $O(n,n|{\bf R})$ on ${\cal T}_n$ by
considering a (non-linear)  embedding of   ${\cal T}_n$ into the Grassmann
algebra  ${\cal F}_n$ with generators $a^1,...,a^n$. To every $\theta
\in{\cal T}_n$ we assign the element 
\begin {equation}
\hat {\theta}=\exp ({1\over 2}a^i\theta _{ij}a^j)\in {\cal F}_n. 
\end {equation}
 Let us introduce operators $\hat {a}^j$and $\hat {b}_k$ transforming
$\omega \in {\cal F}_n$ into $a^j\omega$ and ${\partial \omega \over
\partial a^k}$ respectively. (One can consider ${\cal F}_n$ as a Fock
space; then these operators are creation and annihilation operators.) The
operators $\hat {a}^j\hat {b}_k$ determine an irreducible representation
of the Clifford algebra; in other words, they obey the 
canonical anticommutation
relations: $[\hat {a}^j,\hat {a}^k]_+=[\hat {b}_j,\hat {b}_k]_+=0,\ \ \
[\hat {a}^j,\hat {b}_k]_+=\delta _k^j$. Automorphisms of the Clifford algebra
(linear canonical transformations) constitute a group, isomorphic to
$O(n,n|{\bf R})$. We obtain a projective action of $O(n,n|{\bf R})$ on
${\cal F}_n$ by assigning to every automorphism $\alpha$ an operator
$U_{\alpha}: {\cal F}_n\rightarrow {\cal F}_n$  by the formula $\alpha
(a)\omega =U_{\alpha}(a\omega )U_{\alpha}^{-1}$, where $a$ is an arbitrary
element of the Clifford algebra. Restricting this action to ${\cal T}_n
\subset {\cal F}_n$ we obtain the action (3). Indeed, the element $\hat
{\theta }$ satisfies the equation 
$$\hat {b}_i\hat {\theta}=(\theta _{ij} \hat {a}j)\hat {\theta }  .$$
If 
\begin {equation}
\alpha (\hat {a}^i)=A_j^i\hat {a}^j+B^{ij}\hat {b}_j
\end {equation}
\begin {equation}
\alpha (\hat {b}_i)=C_{ij}\hat {a}^j+D_i^j\hat {b}_j,
\end {equation}
 then $\hat {\theta }^{\prime}=U_{\alpha}\hat {\theta}$ obeys 
\begin {equation}
(c_{ij}\hat {a}^j+D_i^j\hat {b}_j)\hat {\theta}^{\prime}=(\theta
_{ik}A_j^k\hat {a}^j+\theta _{ik}B^{kj}\hat {b}_j)\hat {\theta }^{\prime}.
\end {equation}
In the case when the matrix $D_i^j=\theta _{ik}B^{kj}$ is invertible one
can find $\hat{\theta}^{\prime}$ in the form $\exp ({1\over 2}a^i\theta
_{ij}^{\prime}a^j)$ . The expression for $\theta^{\prime}$ coincides
with (3).

  The description of the action of $O(n,n|{\bf R})$ in terms of the Grassmann
algebra is useful if we would like to relate our theorem to known results
about the $K$-groups and the cyclic homology of noncommutative tori. 

  Let us consider the Grassmann algebra  ${\cal F}_n^{\star}$ dual to the
Grassmann algebra ${\cal F}_n$, and the integral lattice ${\cal
F}_n^{\star}({\bf Z})$  in ${\cal F}_n^{\star}$. The group $K_0(
{A}_{\theta})$ can be identified with the even part $\Lambda ^{{\rm
even}}$ of this lattice and the group $K_1({A}_{\theta}) $ can be
identified with the odd part of it \cite{E}. 
(Recall that a Grassmann algebra has
a natural ${\bf Z}_2$-grading.)  The action of $SO(n,n|{\bf R})$ on ${\cal
F}_n$  induces an action of this group on ${\cal F}_n^{\star}$. The
integral lattice in  ${\cal F}_n^{\star}$ is carried into itself by the
action of  $SO(n,n|{\bf Z})\subset  SO(n,n|{\bf R})$. Therefore we obtain an
action of $SO(n,n|{\bf Z})$ on $K_0$ and on $K_1$. The canonical trace
$\tau$ on ${A}_{\theta}$ determines a group homomorphism $\hat
{\tau}$ of $K_0({A}_{\theta})$ into ${\bf R}$; if an element of
$K_0({A}_{\theta})$ is represented by $x\in \Lambda ^{{\rm even}}$
then $\hat {\tau}(x)$ can be calculated as the scalar product $<\hat
{\theta},x>$ where $\hat {\theta}$ is defined by the formula (5). Taking
into account that $K$-groups of Morita equivalent tori coincide, we obtain
a complete agreement between this expression for $\hat {\tau}$ and our
theorem ( see \cite{S} for more details).

  It is interesting to notice that the embedding of ${\cal T}_n$ into
${\cal F}_n$ described above can be used to define a natural completion of
${\cal T}_n$. Namely, we should consider the closure $\bar {{\cal T}}_n$
of ${\cal T}_n$ embedded in ${\cal F}_n$. The projective action of
$O(n,n|{\bf R})$ on ${\cal F}_n$ induces a projective action of this group
on $\bar {{\cal T}}_n$. (Let us emphasize that there is no ambiguity in
the action of $O(n,n|{\bf R})$ on $ {\cal T}_n$.)

\medskip

{\bf Acknowledgements:}

We would like to thank M. Ratner, G. Prasad, and D. G. James for helpful
comments about $O(n,n|{\bf Z})$. We also thank F. Boca for pointing out an
incorrect assertion in a previous version of this paper, and we thank
G. Elliott for discussion about the difference between the smooth
and C*-algebraic tori.

\noindent
e-mail: rieffel@math.berkeley.edu

\noindent  
e-mail: schwarz@math.ucdavis.edu

                                           \end {document}